\documentstyle[12pt]{article}
\newcommand{\const}{\mathop{\rm const}\limits}

\newcommand{\supp}{\mathop{\rm supp}\limits}

\newcommand{\Div}{\mathop{\rm Div}\limits}

\begin{document}

\begin{center}

{\bf SOLVABILITY OF NAIVER-STOKES EQUATIONS   \\
IN SOME REARRANGEMENT INVARIANT SPACES.} \par

\vspace{4mm}

 $ {\bf E.Ostrovsky^a, \ \ L.Sirota^b } $ \\

\vspace{4mm}

$ ^a $ Corresponding Author. Department of Mathematics and computer science, Bar-Ilan University, 84105, Ramat Gan, Israel.\\
\end{center}
E - mail: \ galo@list.ru \  eugostrovsky@list.ru\\
\begin{center}
$ ^b $  Department of Mathematics and computer science. Bar-Ilan University,
84105, Ramat Gan, Israel.\\

E - mail: \ sirota3@bezeqint.net\\

\vspace{3mm}
                    {\bf Abstract.}\\

 \end{center}

 \vspace{4mm}

 We prove that the multidimensional dimensional initial value problem for the Navier-Stokes equations is globally well-posed in
the so-called Moment and Grand Lebesgue Spaces (GLS), and give some a priory estimations for solution in this spaces. \par

 We consider separately the cases of small initial value solution (local solution)  and global  solution.\par

% We consider separately the cases of short-time  and long-time solutions.\par

\vspace{4mm}

 \par

\vspace{4mm}

  \par

\vspace{4mm}

\section{Notations. Statement of problem.}

\vspace{3mm}

{\bf Statement of problem.} \par

\vspace{3mm}

We  consider in this article the  initial value problem  for the multivariate Navier-Stokes (NS) equations

$$
\partial{u}_t - \Delta u + (u \cdot \nabla)u = \nabla P, \ x \in R^d, \ d \ge 3,   \ t > 0; \eqno(1.1)
$$

$$
\Div (u) = 0, \ x \in R^d, \ t > 0; \eqno(1.2)
$$

$$
u(x,0) = u_0(x), \ x \in R^d. \eqno(1.3)
$$
 Here as ordinary

$$
x = (x_1,x_2,\ldots,x_k,\ldots,x_d) \in R^d ; \  ||x||: = \sqrt{\sum_{j=1}^d x_j^2}.
$$
and $ u = = u(t) = u(t,\cdot)  = u(x,t) $ denotes the (vector) velocity of fluid in the point $ x $ at the time $  t,  \ P $   is
represents the pressure. \par
 Equally:

 $$
\partial{u_i}/\partial t  = \sum_{j=1}^d \partial^2_{x_j} u_i - \sum_{j=1}^d u_j \partial_{x_j} u_i +
+ \partial_{x_i} P,
 $$

$$
\sum_{j=1}^d \partial_{x_j} u_j = 0, \ u(x,0) = u_0(x),
$$

$$
\Div u = \Div \vec{u} = \Div \{ u_1, u_2, \ldots, u_d \} = \sum_{k=1}^d \frac{\partial u_k }{\partial x_k} = 0
$$
in the sense of distributional derivatives.\par
 As long as

 $$
 P = \sum \sum_{j,k = 1}^d R_j R_k (u_j \cdot u_k),
 $$
where $ R_k = R_k^{(d)} $ is the $ k^{th} \ d \- $ dimensional Riesz transform:

$$
 \ R_k^{(d)}[f](x) = c(d) \lim_{\epsilon \to 0+} \int_{ ||y|| > \epsilon} ||y||^{-d} \Omega_k(y) \ f(x-y) \ dy,
$$

$$
c(d) = -\frac{\pi^{(d+1)/2}}{\Gamma \left( \frac{d+1}{2} \right) }, \ \Omega_k(x)= x_k /||x||,
$$
the system  (1.1) - (1.3) may be rewritten as follows:

$$
\partial{u}_t = \Delta u + (u \cdot \nabla)u  +   Q \cdot \nabla \cdot (u \otimes u), \ x \in R^d, \ t > 0; \eqno(1.4)
$$

$$
\Div (u) = 0, \ x \in R^d, \ t > 0; \eqno(1.5)
$$

$$
u(x,0) = u_0(x), \ x \in R^d, \eqno(1.6)
$$
where   $ Q $ is multidimensional  Helmholtz-Weyl projection
operator, i.e., the $ d \times d  $ matrix pseudo-differential operator in $ R^d $ with the matrix symbol

$$
a_{i,j}(\xi) = \delta_{i,j} - \xi_i \xi_j /||\xi||^2, \ \delta_{i,j} = 1, i = j; \delta_{i,j} = 0, \ i \ne j.
$$

\vspace{4mm}
 {\it We will understand henceforth as a capacity of the solution (1.4) - (1.6) the vector - function  } $ u = \vec{u} =
 \{ u_1(x,t),  u_2(x,t), \ldots,  u_d(x,t) \} $ {\it  the so-called mild solution,} see \cite{Muira1}.  \par
\vspace{4mm}

 Namely, the vector- function  $ u = u(t)  $ satisfies almost everywhere in the time $ t $ the following {\it non-linear
 integral - differential equation: }

$$
u(t) = e^{t \Delta} u_0 + \int_0^t e^{(t-s)\Delta } [ (u \cdot \nabla)u(s)  +   Q \cdot \nabla \cdot (u \otimes u)(s) ] ds, \eqno(1.7)
$$
where the operator $  \exp(s \Delta) $ is the classical integral operator with heat kernel: \par

$$
e^{t \Delta}[ u_0](x,t) = (4 \pi t)^{-d/2} \int_{R^d} u_0(y) \ \exp \left( - \frac{|| x - y ||^2}{4 t} \right) \ dy. \eqno(1.8)
$$

\vspace{3mm}

   More results about the existence, uniqueness, numerical methods, and a priory estimates in the different Banach function spaces:
 Lebesgue-Riesz $  L_p, $ Morrey, Besov for this solutions see, e.g. in  \cite{Cui1}- \cite{Temam1}.  The first and besides famous
 result belong to J.Leray  \cite{Leray1}; it is established there in particular the {\it global in time}
 solvability and  uniqueness of NS system in the space $ L_2(R^d)  $ and was obtained  a very interest a priory estimate for solution.\par
  The immediate predecessor for offered article is  the article of Shangbin Cui \cite{Cui1}; in this article was considered the case
  $  u_0 \in L_2(R^d) \cap L_p(R^d), \ p > d.$  See also celebrate works of Y.Giga
  \cite{Giga1} - \cite{Giga4} and T.Kato \cite{Kato1} - \cite{Kato2}.\par

\vspace{4mm}

{\bf  Our purpose in this report is to generalize aforementioned results  about solvability of NS system
on the wide class of rearrangement invariant spaces and to
obtain some useful a priory  estimates for this solution.}  \par

 \vspace{3mm}

 This estimates allow us to  establish some new properties of solution and develop numerical methods. \par

\vspace{3mm}

{\bf Grand Lebesgue Spaces.} \par

\vspace{3mm}

 We recall here briefly the definition and some simple properties  of the so-called
(Bilateral) Grand Lebesgue Spaces (GLS); more detail presentment see, e.g. in \cite{Fiorenza3},
\cite{Kozachenko1},  \cite{Liflyand1},  \cite{Ostrovsky1}, \cite{Ostrovsky2}. \par

 Let $ (X, \Sigma,\mu) $ be a measure space; in the considered problem $  X = R^d $ with Lebesgue measure $ d \mu = dx. $

For $a$ and $b$ constants, $1 \le a < b \le \infty,$ let $\psi =
\psi(p),$ $p \in (a,b),$ be a continuous log-convex positive
function such that $\psi(a + 0)$ and $\psi(b-0)$ exist, with
$\max\{\psi(a+0), \psi(b-0)\} = \infty$ and $\min\{\psi(a+0),
\psi(b-0)\}> 0.$

The (Bilateral) Grand Lebesgue Space (in notation  GLS =  BGLS)
$  G_X(\mu; \psi; a,b) = G_X(\psi;a,b) = G(\psi; a,b) = G(\psi) $
is the space of all measurable functions $ h: X \to R $ endowed with the norm

$$
||h||G(\psi) \stackrel{def}{=}\sup_{p \in (a,b)}||h||_p/\psi(p),
 \quad || h||_p = \left[\int_X |h(x)|^p \ d\mu(x) \right]^{1/p}. \eqno(1.9)
$$

 The  $G(\psi)$ spaces with $ \mu(X) = 1$ appeared in \cite{Kozachenko1};
it was proved that in this case each $ G(\psi) $ space coincides
with certain exponential Orlicz space, up to norm equivalence.
Partial cases of these spaces were intensively studied, in
particular, their associate spaces, fundamental functions
$\phi(G(\psi; a,b);\delta),$ Fourier and singular operators,
conditions for convergence and compactness, reflexivity and
separability, martingales in these spaces, etc.; see, e.g. in
\cite{Fiorenza3}, \cite{Fiorenz4},  \cite{Fiorenz5}, \cite{Kozachenko1},
\cite{Iwaniec1}, \cite{Iwaniec2}, \cite{Jawerth1},  \cite{Liflyand1},  \cite{Ostrovsky1}, \cite{Ostrovsky2}
 etc. \par

 These spaces are also Banach and moreover rearrangement invariant
(r.i.). The BGLS norm estimates, in particular, Orlicz norm estimates for
measurable functions, e.g., for random variables are used in the theory of PDE, probability in Banach spaces,
in the modern non-parametrical statistics, for
example, in the so-called regression problem. \par

We note that the $ G(\psi) $ spaces are also interpolation spaces
(the so-called  $\Sigma$-spaces). However, we hope that our direct
representation of these spaces is of certain convenience in both
theory and applications. A natural question arises what happens if
the spaces other than $L_p$ are used in the definition. Indeed,
this is possible and might be of interest, but, for example, using
Lorenz spaces in this capacity leads to the same object. \par

\vspace{3mm}

{\bf Remark 1.1.} If we define the degenerate $ \psi_r(p), \ r = \const ≥ 1 $ function as
follows:

$$
\psi_r(p) = \infty, \ p \ne r; \  \psi_r(r) = 1, \eqno(1.10)
$$
and agree $  C/\infty = 0, C = \const > 0, $ then the $ G\psi_r(·) $ space coincides with the
classical Lebesgue space $ L_r. $ \par
 Thus, the Grand Lebesgue Spaces are direct generalization of the classical Lebesgue - Riesz spaces. \par

\vspace{3mm}

{\bf Remark 1.2.}  We will denote in the case $  X = R^d  $ by
$  G^0_X(\mu; \psi; a,b) = G^0_X(\psi;a,b) = G^0(\psi; a,b) = G^0(\psi) $
the subspace of the space
$  G_X(\mu; \psi; a,b) = G_X(\psi;a,b) = G(\psi; a,b) = G(\psi) $
which  consists on all the functions $ \{  h = h(x) \}, x \in R^d  $ from this set such that

$$
\Div h = 0. \eqno(1.11)
$$
 Analogously, the space $  L_p^0 $ consists on all the functions from the space $  L_p $ with zero  divergence. \par

 \vspace{3mm}

 {\it In what follows we will suppose } $ \Div u_0 = 0; $  {\it  therefore, }   $ \Div u_(x,t) = 0 $  {\it for all the
 values $  t  $  for  which the solution $  u(x,t)   $  there exists.} \par

\vspace{3mm}

{\bf Remark 1.3.} Multidimensional case. \par
 Let $ u = \vec {u} = \{ u_1(x), u_2(x), \ldots, u_d(x)  \}  $ be measurable vector - function. We can define
as ordinary the $ G\psi $  norm of the function $ u $ by the following way:

$$
||u||G\psi := \max_{k=1,2, \ldots,d} || u_k|| G\psi.
$$

\vspace{3mm}

{\bf Remark 1.4.} Natural choice. \par
 Let $ v = v(x) \ne 0, \ x \in X  $ be some (measurable) function for which there exist two constants
 $ a, b: \ 1 \le a < b \le \infty $  such that

 $$
 \forall p \in (a,b) \ \Rightarrow || v ||_p < \infty.
 $$
 The $ \psi = \psi_{(v) }(p)  $ function of a view

$$
\psi_{(v) }(p) = || v ||_p, \ p \in (a,b)
$$
 is said to be {\it  natural function } for the function $  v(\cdot). $\par

\vspace{3mm}

 Analogously, for the {\it family} of a functions $ v = v(x, \alpha), \ \alpha \in A, \ A $ is an arbitrary set, the
 {\it  natural function } $ \psi = \psi_{(v) }(p) $ may be defined as follows:

$$
\psi_{(v) }(p) =  \sup_{\alpha \in A} || v (\cdot, \alpha)||_p, \ p \in (a,b),
$$
 if there exists.\par

 Obviously,

$$
\sup_{\alpha \in A} ||v(\cdot,\alpha)|| G\psi_{(v) } = 1.
$$

 This approach is very convenient, e.g., in the theory of Probability and Statistics. \par

\vspace{4mm}

{\bf Moment Rearrangement Invariant Spaces.} \par

\vspace{3mm}

  Let $ (X = R^d, \ ||\cdot||X) $ be any r.i. space, where $ X $ is linear subset on the space of all
 measurable function $ R^d \to R $ over our measurable space
 $ (T,M,\mu) $ with norm $ ||\cdot||X. $ Recall the following definition,
 see, e.g. \cite{Ostrovsky107}, \ \cite{Ostrovsky108}, where are described some applications of these
 spaces in the Approximation Theory and in the Theory of Partial Differential Equations. \par

\vspace{4mm}

 We will say that the space $ X $ with the norm $ ||\cdot||X $ is {\it moment
 rearrangement invariant space,} briefly: m.r.i. space, or
$ X =(X, \ ||\cdot||X) \in m.r.i., $
 if there exist a real constants $ a, b; 1 \le a < b \le \infty, $ and some {\it rearrangement invariant norm }
  $ < \cdot > $ defined on the space of a real functions defined on the interval $ (a,b), $ non necessary to be
  finite on all the functions, such that

  $$
  \forall f \in X \ \Rightarrow || f ||X = < \ h(\cdot) \ >, \ h(p) = |f|_p. \eqno(1.12)
  $$

   We will say that the space $ X $ with the norm $ ||\cdot||X $ is {\it weak moment
 rearrangement space,} briefly, w.m.r.i. space, or $ X =(X, \ ||\cdot||X) \in w.m.r.i.,$
 if there exist a constants $ a, b; 1 \le a < b \le \infty, $ and some {\it functional } $ F, $ defined on
 the space of a real functions defined on the
interval $ (a,b), $ non necessary to be finite on all the functions, such that

  $$
  \forall f \in X \ \Rightarrow || f ||X = F( \ h(\cdot) \ ), \ h(p) = |f|_p.
  $$

    We will write for considered w.m.r.i. and m.r.i. spaces $ (X, \  ||\cdot||X) $

   $$
       (a, b) \stackrel{def}{=} \supp(X),
   $$
   (“moment support”; not necessary to be uniquely defined)
   and define for other such a space $ Y = (Y, \ ||\cdot||Y ) $ with
   $ (c,d) = \supp(Y) $
   $$
   \supp(X) >> \supp(Y),
   $$
    iff $ \min(a,b) > \max(c,d). $ \par
     It is obvious that arbitrary m.r.i. space is r.i. space.\par

\vspace{3mm}

\section{ Some Notations, with Clarification.} \par

\vspace{3mm}

 As ordinary, for the measurable function $ x \to u(x), \ x \in R^d $

 $$
||u||_p = \left[ \int_{R^d} |u(x)|^p \ dx  \right]^{1/p}. \eqno(2.1)
 $$
 The so-called {\it mixed,  } or equally {\it anisotropic } $ (p_1, p_2) $   norm  $ ||u||^*_{p_1,p_2} $ for the function of "two"
variables $ u = u(x,t), \ x \in R^d, \ t \in R^1_+ $  is defined as follows:

$$
||u||^*_{p_1,p_2} = \left( \int_{R^d} \left[ \int_0^{\infty} |u(x,t)|^{p_1} \ dt  \right]^{p_2/p_1} \ dx \right)^{1/p_2}. \eqno(2.2)
$$

 The correspondent {\it mixed,  } or  {\it anisotropic } Grand Lebesgue (Lebesgue - Riesz) spaces was introduced in
\cite{Ostrovsky110}.   For the positive function of two variables $ \theta = \theta(p,r) $  defined on the set $  D $
the norm of a function in this space is defined by formula

$$
||u||^* G \theta = \sup_{(p,r) \in D  } \left[ \frac{||u||^*_{p,r}}{\theta(p,r)} \right]. \eqno(2.3)
$$

\vspace{5mm}

  The following functional $ u \to \kappa_p^{(d)}(u) = \kappa_p(u)  $ in the case  $ d=3 $  was introduced
 and used by Shangbin Cui in \cite{Cui1}:

$$
 \kappa_p(u) = \kappa_p^{(d)}(u) := ||u||_p^{\frac{p(d-2)}{d(p-2)}} \ ||u||_2^{ \frac{2(p-d)}{d(p-2)} }. \eqno(2.4)
$$
 This functional is scaling-dilation invariant. Indeed, denote

 $$
 T_{\lambda} [u](x) = \lambda \ u(\lambda x), \ \lambda \in R, \ \lambda \ne 0,
 $$
then

$$
\kappa_p^{(d)}(T_{\lambda} [u]) =  \kappa_p^{(d)}(u).
$$

\vspace{5mm}

 Define also

 $$
W = W_{d,p} = W_{d,p}(u) := \int_{R^d} |u(x)|^{p-2} \ |\nabla u|^2 \ dx;
 $$

$$
K_S(d,p) := \pi^{-1/2} \ d^{-1/p} \ \left( \frac{p-1}{d-p} \right)^{ (p-1)/p } \
\left\{ \frac{\Gamma(1+d/2) \ \Gamma(d)}{\Gamma(d/p) \ \Gamma(1+d - d/p) } \right\}^{1/d}. \eqno(2.5)
$$

 The function $ K_S(d,p) $ is the optimal (i.e. minimal) value in the famous Sobolev's inequality

 $$
 ||\phi||_r \le K_S(d,q) \ || \nabla \phi  ||_q, \ 1 \le q < d, \ \frac{1}{r} = \frac{1}{q} -  \frac{1}{d}, \ r \ge 1, \eqno(2.6)
 $$
see  Bliss \cite{Bliss1}, (1930); Talenti, \cite{Talenti1}, (1995).\par

\vspace{5mm}

 Further, denote

$$
A_{d,p}:= \left( \frac{p+d}{p}  \right)^{\frac{p+d}{(p-d)}}, \hspace{5mm} \ B_{2.1}(d,p) := K_S^2(d,2d/3) \ p^2/4. \eqno(2.7)
$$

 Note that

 $$
 K_S(d,2d/3) = 2^{1/d} \cdot \pi^{-\frac{d+1}{2d} } \cdot (2-3/d) \cdot (2d-3)^{ -3/2d } \times
 $$

$$
\left\{ \frac{\Gamma(1+d/2) \ \Gamma(d)}{\Gamma(d-1/2)} \right\}^{1/d}, \ d = 3,4,5, \ldots. \eqno(2.8)
$$
 Since the number $  d  $ is integer, the expression for $ K_S(d,2d/3) $ may be calculated in explicit view.
 For instance,
 $$
 \Gamma(d) = (d-1)!, \ \Gamma(d-1/2) = \frac{\sqrt{\pi}}{2^{d-1}}  \ (2d - 3)!!.
 $$
 So, if the dimension $ d $ is even number, then

$$
 K_S(d,2d/3) = 2^{1/d} \cdot \pi^{-\frac{d+1}{2d} } \cdot (2-3/d) \cdot (2d-3)^{ -3/2d } \times
$$

$$
\left\{ \frac{2^{d-1} \ (d/2)! \  (d-1)!}{\sqrt{\pi} \ (2d-3)!!}  \right\}^{1/d}.
$$
 In the opposite case, i.e. when $  d  $ is odd number,

$$
 K_S(d,2d/3) = 2^{1/d} \cdot \pi^{-\frac{d+1}{2d} } \cdot (2-3/d) \cdot (2d-3)^{ -3/2d } \times
$$

 $$
 \left\{ \frac{2^{(d-3)/2} \ d!! \ (d-1)! }{(2d-3)!!}  \right\}^{1/d}.
 $$
  The behavior  of the variable $ K_S(d,2d/3)  $ as $ d \to \infty  $ may be obtained from the Stirling's formula:

 $$
 K_S(d,2d/3) \sim \sqrt{ \frac{d}{2 \ \pi \ e}}.
 $$

  For example, at $ d=3 \ $ (the most important case in practice)

$$
K_S(3,2)= \frac{1}{3} \cdot \sqrt[3] {\frac{2}{\pi^2} }.
$$

\vspace{4mm}

 We will use the following elementary inequality

 $$
 vw \le A(d,p) \ v^{^{ \frac{2p}{p-d} } } + 0.5 \ w^{ \frac{2p}{p+d }}, \eqno(2.9)
 $$
where $ p > d,  \ v,w > 0. $  Therefore,

$$
||u||_p^{1+ (p-d)/2  } \cdot \left( \int_{R^d} |u(x)|^{p-2} \ |\nabla u(x)|^2 \ dx    \right)^{\frac{p+d}{2p}} \le
$$

$$
A(d,p) || u ||_p^{ \frac{p(p-d+2)}{p-d}} + \frac{1}{2} \ \left( \int_{R^d} |u(x)|^{p-2} \ |\nabla u(x)|^2 \ dx    \right). \eqno(2.10)
$$

\vspace{3mm}

$$
\tilde{\omega}(d) := \frac{4 \pi^{d/2-1}}{\Gamma(d/2)}. \hspace{4mm}I(p) :=  \frac{1}{2\sqrt{\pi}} \
\Gamma \left( \frac{1}{2} -  \frac{1}{2 p} \right)  \ \Gamma \left( \frac{1}{2 p}  \right). \eqno(2.11)
$$

\vspace{3mm}

$$
c(d) = -\frac{\pi^{(d+1)/2}}{\Gamma \left( \frac{d+1}{2} \right) }. \hspace{6mm} \Omega_k(x)= x_k  / ||x||. \eqno(2.12)
$$

\vspace{3mm}

$$
x = (x_1,x_2,\ldots,x_k,\ldots,x_d) \in R^d  \ \Rightarrow ||x|| = \sqrt{\sum_{j=1}^d x_j^2}.
$$

$$
K_R(d,p) =  c(d) \ \cdot \frac{p}{p-1} \cdot  \tilde{\omega}(d) \cdot I(p), \ p > 1. \eqno(2.13)
$$

\vspace{3mm}

 The explicit view for Riesz's transform has a view
$$
 R_k[f](x)= R_k^{(d)}[f](x) = c(d) \lim_{\epsilon \to 0+} \int_{ ||y|| > \epsilon} ||y||^{-d} \Omega_k(y) \ f(x-y) \ dy.
$$

 It is known, see \cite{Okikiolu1}, p. 415-418, that $ || R_k||(L_p \to L_p) \le  $

$$
 c(d) \ \cdot \frac{p}{p-1} \cdot \int_{\Sigma(d)} |x_1| d \sigma_d \cdot \int_0^{\infty}t^{-1/p} (1+t^2)^{-1/2} dt =
 K_R(d,p), \ p > 1.
$$

 Here $ \Sigma(d) $ is an unit sphere in the space $ R^d $ and $ d \sigma_d $ is an element of its area. \par
Note that the last estimate is not improvable even in the case $  d=1, $  where the Riesz transform coincides with
Hilbert transform, for which the norm estimates $ (L_p \to L_p) $  is computed by S.K.Pichorides \cite{Pichorides1}. \par

 Ultimate result in this direction belongs to T.Iwaniec  and  G.Martin  \cite{Iwaniec3}:  the value
$ ||R_k||(L_p \to L_p) $ does not dependent on the dimension $ d $ and coincides with the Pichorides constant:

$$
||R_k||(L_p \to L_p) = \cot \left( \frac{\pi}{2p^*}  \right), \ p^* = \max(p, \ p/(p-1) ), \ p > 1. \eqno(2.14)
$$
 T.Iwaniec  and  G.Martin considered  also the vectorial Riesz transform. \par
 See for additional information  \cite{Banuelos1},  \cite{Stein1},  chapter 2, section 4; \cite{Taylor1}, chapter 3. \par

 Put also

 $$
 C_{2.7}(d,p) := 4^{-1} \ p^2 \ K^2_S(d,2d/3) \  K^2_R(d,p) \ (d^2 + d),
 $$

$$
 C_{7.7}(d,p) = A(d,p) \cdot \left[ C_{2.7}(d,p) \right]^{ \frac{2p}{p-d}}, \eqno(2.15)
$$
we borrow notations from \cite{Cui1} after estimation and specification.\\

\vspace{4mm}

 Introduce also the following function:
\vspace{3mm}

$$
Z = Z_{a,b}(x, y; p) := x^{ \frac{a(b-p)}{p(b-a)} } \cdot y^{ \frac{b(p-a)}{p(b-a)} }, \eqno(2.16)
$$

 $  1 < a < b < \infty; \ p > 1, \ x,y \in (0, \infty). $ \par

  This function has a following  sense: if

$$
 1 < a < b < \infty;  \ f \in L_a \cap L_b, \ p \in (a,b),
$$
then

 $$
 ||f||_p \le Z_{a,b}(||f||_a, ||f||_b; p). \eqno(2.17)
 $$
 The last inequality may be deduced from the H\"older's inequality. \par

\vspace{4mm}

\section{ Solution for small initial data }

\vspace{4mm}

{\it We suppose in this section that the initial function $ u_0 = u_0(x) $ belong to some Grand Lebesgue Space
 $  G\psi $  such that } $  \Div u_0 = 0 $ {\it and } $ d \in \supp \psi.  $ \par

\vspace{3mm}

 But we do not assume here that $ 2 \in \supp;  $ this case will be considered further.\par

\vspace{3mm}

 It is known in the case when $ u_0 \in L_r, \ r \ge d   $  and when the initial norm $ ||u_0||_r  $
 is sufficiently small, then the NS equation has a unique global (smooth) solution;
 see for example  \cite{Cui1}, \cite{Kato1} - \cite{Kato2}. \par

\vspace{3mm}

 We generalize these results on the Grand Lebesgue spaces, calculating passing the constants values. \par

\vspace{3mm}

 Let us denote

 $$
J = \left\{ p: ||u_0||_d < \frac{1}{2C_{7.7}(d,p)}  \right\} \eqno(3.1)
 $$
and introduce the following $ \psi \ - $ function:

$$
\tilde{\psi}(p) = \psi(p), \  p \in \supp \psi \cap J; \  \tilde{\psi}(p) = 0 \eqno(3.2)
$$
otherwise.

\vspace{4mm}

{\bf Theorem 3.1.} Suppose $ \Div u_0 = 0, \ d \in \supp \psi  $  and $ \supp \tilde{\psi} \ne \emptyset. $ Then
the global in time solution of NS system $ u(t) $ there exists with monotonically decreased norm
$ ||u(t)||G^0 \tilde{\psi}  $ and moreover

$$
\sup_{t \ge 0} ||u(t)||G^0 \tilde{\psi} = ||u_0|| G\psi. \eqno(3.3)
$$

\vspace{4mm}

{\bf Proof.} We follow  Shangbin Cui \cite{Cui1} specifying passing the "constants" values but omitting some
hard calculations. \par

\vspace{3mm}

{\bf 1.} If

$$
||u_0||_d \le \frac{1}{2C_{2.7}(d,d)}, \eqno(3.4)
$$
then there exists and is unique the global in time  solution of NS system $ u(t), \ t > 0 $ such that
the function $ t \to ||u(t)||_d  $ is monotonically decreasing. \par

\vspace{3mm}

{\bf 2.} In what follows in this section we suppose $ p \in \supp \psi. $ We conclude using the estimates  for
Riesz transform:

\vspace{3mm}

$$
\frac{1}{p} \ \frac{d}{dt} ||u||_p^p + \int_{R^d} |u(x,t)|^{p-2} \ |\nabla u(x,t)|^2 \ dx \le
$$

$$
 C_{2.7}(d,p) \cdot ||u||_p^{ 1+ (p-d)/2 } \cdot \left( \int_{R^d} |u(x,t)|^{p-2} \cdot |\nabla u(x,t)|^2 \right)^{\frac{p+d}{2p}}.
 \eqno(3.5)
$$

\vspace{3mm}

{\bf 3.} We obtain after Shangbin Cui  \cite{Cui1} by means of constant computation and using the expression for $ A(d,p):  $

$$
\frac{1}{p} \ \frac{d}{dt} ||u||_p^p + 0.5 \ \int_{R^d} |u(x,t)|^{p-2} \ |\nabla u(x,t)|^2 \ dx \le
$$

$$
C_{7.7}(d,p) ||u||_p^{ \frac{p(p-d+2)}{p-d}}. \eqno(3.6)
$$

\vspace{3mm}

{\bf 4.} Therefore, if

 $$
 ||u_0||_d < \frac{1}{2C_{7.7}(d,p)},
 $$
then the function

$$
t \to ||u(t)||_p, \ t \ge 0, \ p \in J,
$$
 is monotonically decreasing, which is equivalent to the assertion of theorem 3.1. \par

\vspace{4mm}

\section{Global solution.}

\vspace{3mm}

 It is proved in \cite{Cui1}  that if
 $$
  u_0(\cdot) \in L_2^0 \cap L_b^0  \eqno(4.1)
 $$
 for  some $ b  \ge d,  $ then  the  Navier - Stokes equations, more precisely, the system of  Navier - Stokes equations
 (1.1) - (1.3) has unique global in time smooth solution $ u = u(x) = u(x,t). $ But  if the condition (4.1) is satisfied,
 we still have to take admit that the function $ u_0(\cdot) $ belongs to some Grand Lebesgue space.  In detail,
  there exists a function $ \psi = \psi_b(p)  $ with support $ [2, b) $ such that  $ \forall t \ge 0 \ u(\cdot) \in G\psi.  $ \par

  Indeed, let $ y_2:=||u||_2  < \infty  $ and  $ y_b:=||u||_b  < \infty;  $ it follows from the H\"older's inequality that
 for any value $ p \in (2,b) $

 $$
 ||u||_p \le Z_{2,b}(y_2, y_b; p) =:\psi_b(p).
 $$

{\it  Therefore, it is  more than natural to suppose in this section  that
the initial value function  $  u_0 = u_0(x)  $  belongs
to some  Grand Lebesgue Space } $ G\psi, \ \supp \psi = [2,b) $ {\it with } $ b > d. $\par

% "natural function"

\vspace{3mm}

{\bf Theorem 4.1.} Let the initial value function  $  u_0 = u_0(x)  $  belongs
to some  Grand Lebesgue Space  $ G\psi, $  such that $ \supp \psi = [2,b) $ {\it with } $ b > d. $
 Define a new function

 $$
 \psi_{(\kappa)}(p) := \psi(p) \cdot \max(1, \kappa_p^{2 d/p}(u_0)). \eqno(4.2)
 $$

 Proposition:

 $$
 \sup_t || u(t) || G \psi_{(\kappa)}  \le 1.\eqno(4.3)
 $$

\vspace{3mm}

{\bf Proof.} Let $ u_0 \in G \psi,  $ then $ ||u_0||_p < \infty   $ and moreover

$$
||u_0||_p \le  \psi(p), \ 2 \le p < b.
$$

  It may be deduced after some (omitted here) calculations based on the  the article of Shangbin Cui  \cite{Cui1} that

$$
  \sup_t ||u||_p \le \max(1, \kappa_p^{2 d/p}(u_0) ) \cdot ||u_0||_p  \le \max(1,\kappa_p^{2 d/p}(u_0)) \cdot \psi(p),
$$
or equally

$$
 \sup_t || u(t) || G \psi_{(\kappa)} =
\sup_t \sup_{ 2 \le p < b} \frac{||u(\cdot,t)||_p}{\psi_{(\kappa)}(p)} \le 1. \eqno(4.4)
$$
 This completes the proof of theorem 4.1.\par

\vspace{4mm}

\section{Solvability in moment rearrangement spaces.}

\vspace{3mm}

{\it  We retain the notations and assumptions of previous section.} \par

\vspace{3mm}

  Let $ (X, \ ||| \cdot ||| X) $ be any moment rearrangement invariant space over $  R^d $
 constructed by means  of auxiliary space $  (V, <  \cdot >) $  with condition  $ \supp V = \supp \psi. $
  Denote

  $$
  h_0(p) = \max(1,\kappa_p^{2d/p}) \cdot \psi(p), \ p \in \supp \psi. \eqno(5.1)
  $$

  \vspace{3mm}

  {\bf Theorem 5.1.}

  $$
  \sup_{ t \ge 0} ||| u(t)|||X \le \hspace{3mm} <  h_0(\cdot) >. \eqno(5.2)
  $$

\vspace{3mm}
{\bf Proof }  is very simple. It follows from theorem 4.1

$$
 ||u(t) ||_p \le \max(1,\kappa_p^{2d/p}) \cdot \psi(p) = h_0(p), \ t \ge 0. \eqno(5.3)
$$

 Since the norm $ < \cdot >  $ is also rearrangement invariant, we deduce taking the norm   $ < \cdot >  $  on both the sides
of inequality (5.3):

$$
< ||u||_p  > \hspace{3mm} \le \ \hspace{3mm} < h_0(\cdot)  >,  \eqno(5.4)
$$
and this estimate is uniform on the variable $  t; t \ge 0. $ \par
 Taking the supremum over variable $ t, $  we get to the  proposition (5.2) of the regarded theorem. \par

\vspace{4mm}

\section{Mixed norm estimates for  solution.}

\vspace{3mm}

 It is known, see \cite{Cui1}, \cite{Kato1}-  \cite{Kato2} that  the global in time solution $ u(x,t) = u(t)  $ obeys the
property

$$
\lim_{t \to \infty} ||u(t)||_p = 0.
$$

 We want in this section to characterize this  feature on the language of anisotropic Grand Lebesgue spaces.\par

 We suppose in this section that $ u_0 \in L_2^0 \cap L_b^0, \ b > d  $ or equally that the initial function $ u_0   $
belongs to some $ G\psi  $ space with  $ \supp \psi = [2, b), \ d < b \le \infty.  $ \par

 A new notations:

 $$
 r = r(p) = r_d(p): = \frac{p(p-d+2)}{p-d}, \ p \in (d,b); \ D = \{ p, \ r(p)  \}, \eqno(6.1)
 $$

$$
\theta_{d, \psi}(p) :=
\left[ \frac{B_{2.1}(d,p)}{p} \right]^{1/r(p)} \cdot \psi^{ \frac{2}{p-d+2}}(d) \cdot
\psi^{ \frac{p-d}{p-d+2}}(p). \eqno(6.2)
$$

\vspace{3mm}

{\bf Theorem 6.1.} Assume $ u_0 \in G \psi;  $ then

$$
|| u||^* G \theta_{d, \psi} \le 1. \eqno(6.3)
$$

\vspace{3mm}

{\bf Proof } is at the same as before: we estimate the norm  $ ||u||^* G \theta_{d, \psi}   $ specifying passing
the "constants"  from the article of Shangbin Cui  \cite{Cui1}.  In detail, let $ p \in \supp \psi;  $ then

$$
\int_0^{\infty} ||u(t)||^{r(p)}_p \ dt \le \frac{B_{2.1}(d,p)}{p} \ ||u_0||_d^{ \frac{2p}{p-d}  } \cdot ||u_0||^p_p. \eqno(6.4)
$$

 Since $ u_0 \in G \psi,  $

$$
||u_0||_d \le \psi(d), \hspace{5mm} ||u_0||_p \le \psi(p). \eqno(6.5)
$$
 We get substituting into (6.4) and taking the root of degree $  r = r(p)  $

$$
||u||^*_{p, r(p)} \le \left[ \frac{B_{2.1}(d,p)}{p} \right]^{1/r(p)} \cdot \psi^{ \frac{2}{p-d+2}}(d) \cdot
\psi^{ \frac{p-d}{p-d+2}}(p) = \theta_{d, \psi}(p).
$$
 It remains to divide  into $ \theta_{d, \psi}(p) $ and take supremum over $ (p, r)\in D. $ \par

 \vspace{4mm}

 \section{Concluding remarks. }

 \vspace{4mm}

{\bf 1.  First example.  }\par
\vspace{3mm}
 We retain here the assumptions and notations of fourth sections. Suppose in addition that

 $$
  \sup_{2 \le p < b} \psi^{1/p} (p) < \infty. \eqno(7.1)
 $$
 This condition is satisfied if for example $  b  = \infty   $ and $ \psi(p) = p^m \ M(p), $ where $ m = \const < \infty $
and $ M(p) $ is positive continuous slowly varying as $ p \to \infty $ function. \par
 We have under condition (7.1)

$$
\psi_{(\kappa)}(p) \asymp \psi(p),
$$
therefore, the proposition of theorem 4.1 may be rewritten in the considered case  as follows:

 $$
 \sup_t || u(t) || G \psi  < \infty. \eqno(7.2).
 $$

\vspace{4mm}

{\bf 2.  Second example.  }\par

\vspace{3mm}

We consider here the mixed norm estimates for  solution, see sixth section. Note that

$$
\theta_{d, \psi}(p) \asymp \psi(p),
$$
even without the condition (7.1).
Following, the proposition of theorem 6.1  has a view
 $$
|| u||^* G \psi < \infty. \eqno(7.3)
$$

\vspace{4mm}

{\bf 3.}  It is known \cite{Giga1}, \cite{Giga2},  \cite{Kato1}, \cite{Kato2}  etc. that in general case, i.e.
when the value $ \epsilon = ||u_0||_d $ is not sufficiently small, then the lifespan of solution of NS equation $ T $
may be finite (short-time solution).  Perhaps, it is  self-contained  interest a quantitative estimate of the value $ T. $ \par
  For the non-linear Schr\"odinger's equation the estimate

 $$
 T \ge  \exp (C /\epsilon)
 $$
was obtained in the recent article \cite{Ikeda1}.\par

\vspace{4mm}
{\bf 4.} At the same considerations may be provided for the NS equations with external force $  f = f(x,t): $

$$
\partial{u}_t = \Delta u + (u \cdot \nabla)u  +   Q \cdot \nabla \cdot (u \otimes u) + f(x,t), \ x \in R^d, \ t > 0;
$$

$$
u(x,0) = u_0(x), \ x \in R^d.
$$
see  \cite{Giga1} -  \cite{Giga4}, \cite{Koch1},  \cite{Masuda1},  \cite{Temam1}.\par

\vspace{4mm}

 {\bf 5.} Analogously to the content of this report may be considered a more general case of abstract (linear or not linear) parabolic
 equation of a view

 $$
 \partial{u}_t  = A u + F(u, \nabla u; x,t)  + f(x,t), \hspace{5mm} u(x,0) = a(x).
 $$
 The detail investigation of this case when the initial condition and external force belong to some
 Sobolev's space may be found, e.g. in \cite{Taylor1} -  \cite{Taylor3}, \cite{Iwashita1}, \cite{Solonnikov1}. \par


\vspace{4mm}

\begin{thebibliography}{99}

\vspace{3mm}

\bibitem{Cui1}
{\sc Shangbin Cui.} {\it Global well-posedness of the 3-dimensional Navier-Stokes initial value problem
in $ L(p) \cap L(2) $ with } $ 3 < p < \infty. $
arXiv:1204.5040v1 [math.AP] 23 Apr 2012

\bibitem{Barraza1}
{\sc Barraza O.} {\it Self-similar solutions in weak $ L_p - $ spaces of the Navier-Stokes equations. } Revista
Mat. Iberoamer., 12(1996), 411 \ – \ 439.

\bibitem{Caffarelli1}
{\sc Caffarelli L., Kohn R. and Nirenberg L.} {\it Partial regularity of suitable weak solutions of the
 Navier-Stokes equations. } Comm. Pure Appl. Math., {\bf 35,} (1082), 771 \ - \ 831.

\bibitem{Calderon1}
{\sc Calder´on C.} {\it Existence of weak solutions for the Navier-Stokes equations with initial data
in $ L(p). $ } Trans. A.M.S., 318(1990), 179 \ – \ 207.

\bibitem{Escauriaza1}
{\sc Escauriaza L., Seregin G., and  Sver´ak V. }$ L(3,\infty) $ {\it -Solutions to the Navier-Stokes Equations
and Backward Uniqueness. } Uspekhi Mat. Nauk, 58( 2003), no.2, 3 \ – \ 44.

\bibitem{Fabes1}
{\sc Fabes E., Johns B. and  Riviere N. } {\it The initial value problem for the Navier-Stokes equations
with data in $ L(p). $ } Arch. Rat. Mech. Anal., 45(1972), 222 \ – \ 240.

\bibitem{Foias1}
{\sc Foias C., Guillope C. and  Temam R.} {\it New a priory estimates for Navier-Stokes equations in dimension 3. }
Comm. in Part. Dif. Eq., {\bf 6,} (1981), 329 \ - \ 359.

\bibitem{Fujita1}
{\sc  Fujita H. and Kato T. } {\it On the Navier-Stokes initial value problem I. } Arch. Ration. Mech.
Anal., 16(1964), 269 \ – \ 315.

\bibitem{Germain1}
{\sc Germain P. } {\it Multipliers, paramultipliers, and weak-strong uniqueness for the Navier-Stokes
equations. } J. Diff. Equations., 226(2006), 373 \ – \ 428.

\bibitem{Giga1}
{\sc Giga  Y. } {\it Solutions of semilinear parabolic equations in $ L^p $ and regularity of weak solutions
of the Navier-Stokes system. } J. Diff. Equations, 62(1986), 186\ – \ 212.

\bibitem{Giga2}
{\sc  Giga Y. and  Miyakawa T. } {\it Navier-Stokes flows in $ R^3 $ with measurea as initial vorticity and
the Morrey spaces. } Comm. P. D. E., 14(1989), 577 \ – \ 618.

\bibitem{Giga3}
{\sc  Giga Y. and  Sohr H. }  {\it Abstract $ L^p \ - $ estimates for the Cauchy  problem with Applications to the
Navier-Stokes equations in  exteroir domains. } Hokkaido University, Preprint, Series 60 on Mathematics, November 1989.

\bibitem{Giga4}
{\sc  Giga Y. and  Sohr H. }  {\it Abstract $ L^p \ - $ estimates for the Cauchy  problem with Applications to the
Navier-Stokes equations in  exteroir domains. } J. Funk. Anal., {\bf 102} (1991), 72 \ -  \ 94.

\bibitem{Iwashita1}
{\sc Iwashita H.} {\it $ L^q - L^r  $  estimates for solution of non-stationary Stokes  equations in exterior  domain
and the Navier-Skokes initial value problems in $ L_q $ spaces. } Math. Ann., {\bf 285}, (1989), 265 \ - \ 288.

\bibitem{Kato1}
{\sc  Kato T. }  {\it Strong $ L_p $ solutions of the Navier-Stokes equations in $ R^m $ with applications to
weak solutions. } Math. Z., 187(1984), 471 \ – \ 480.

\bibitem{Kato2}
{\sc Kato T. and  Ponce G.} {\it Commutator estimates and the Euler and Navier-Stokes equations. }
Comm. P. D. E., 41(1988), 891 \ - \ 907.

\bibitem{Kenig1}
{\sc  Kenig C.E. and Koch G.S. } {\it An alternative approach to regularity for the Navier-Stokes
equations in a critical space. } arXiv:0908.3349.

\bibitem{Koch1}
{\sc Koch H. and Tataru D. } {\it Well-posedness for the Navier-Stokes equations}. Adv. in Math.,
157(2001), 22 \ – \ 35.

\bibitem{Kozono1}
{\sc Kozono H. and Taniuchi Y.} {\it Bilinear estimates in BMO and the Navier-Stokes equations.}
Math. Z., 235(2000), 173 \ – \ 194.

\bibitem{Lemari1}
{\sc Lemari´e-Rieusset P.G. } {\it Weak infinite-energy solutions for the Navier-Stokes equations in
$ R^3,$ } Preprint, 1998.

\bibitem{Lemari2}
{\sc  Lemari´e-Rieusset P.G.} {\it Recent developments in the Navier-Stokes problems.} Research Notes
in Mathematics, Chapman, Hall/CRC, 2002.

\bibitem{Leray1}
{\sc Leray J. } {\it Sur le mouvement d’un liquide visqueux emplissant l’espace. } Acta Math., 63(1934),
193 \ – \ 248.

\bibitem{Masuda1}
{\sc Masuda K.} {\it Weak solutions of Navier-Stokes equations.} T'ohoku Math. J., 36(1984), 623 \ – \ 646.

\bibitem{Muira1}
{\sc  Miura H.} {\it Remarks on uniqueness of mild solutions to the Navier-Stokes equations.} J. Funct.
Anal., 218(2005), 110–129.

\bibitem{Planchon1}
{\sc Planchon F.} {\it Global strong solutions in Sobolev or Lebesgue spaces to the imcompressible
Navier-Stokes equations in} $ R^3,$  Ann. Inst. H. Poincare Anal. Non Lineaire, 13(1996), 319 \ – \ 336.

\bibitem{Seregin1}
{\sc Seregin G. } {\it A certain necessary condition of potential blow up for Navier-Stokes equations.}
arXiv:1104.3615.

\bibitem{Serrin1}
{\sc  Serrin J.} {\it The initial value problem for the Navier-Stokes equations. } In: R.E. Langer, (Ed.),
Nonlinear Problems, 1963, University of Wisconsin Press, Madison, 1963, pp. 69 \ – \ 98.

\bibitem{Solonnikov1}
{\sc Solonnikov V.A.} {\it Estimates for Solutions of non-stationaty Navier - Stokes equations. }
J. Soviet Math., {\bf 8}, (1977),  467 \ - \ 523.

\bibitem{Stein1}
{\sc Stein E. M.} {\it Singular Integrals and Differentiability Properties of Functions.} Princeton,
University Press, (1970),  Princeton, New Jersey.

\bibitem{Taylor1}
{\sc Taylor M.E.} {\it Pseudodifferential Operators.} Princeton, University Press; Princeton, New Jersey,
(1981)

\bibitem{Taylor2}
{\sc Taylor M.E.} {\it Partial Differential Equations I. Linear Equations.} Applied Math. Sciencies, {\bf 117,}
Springer, (1996).

\bibitem{Taylor3}
{\sc Taylor M.E.} {\it Partial Differential Equations III. Non-linear Rquations.} Applied Math. Sciencies, {\bf 117,}
Springer, (1996).

\bibitem{Temam1}
{\sc Temam R.} {\it Navier - Stokes Equations. Theory and Numerical Analysis. } North-Holland Publishing Company.
Amsterdam, New York,Oxford, (1977).

\vspace{9mm}

\bibitem{Banuelos1}
{\sc Ba\"nuelos R. and Osekowski A.} {\it Sharp martingale inequalities and applicationa to Riesz
transforms on manifolds, Lie group and Gauss space. }
arXiv:1305.1492v1 [math.PR] 7 May 2013

\bibitem{Bliss1}
{\sc Bliss G.} {\it  An integral inequality.  } J. London Math. Soc., (1930), vol. 5, 40 \ - \ 46.

\bibitem{Brascamp1}
{\sc Brascamp H.J. and E.H. Lieb E.H.} {\it Best constants in Young’s inequality, its converse and its generalization
to more than three functions.} Journ. Funct. Anal., 20(1976), 151 \ – \ 173.

\bibitem{Fiorenza3}
 {\sc Fiorenza A., and Karadzhov G.E.} {\it Grand and small Lebesgue spaces and
 their analogs.} Consiglio Nationale Delle Ricerche, Instituto per le
 Applicazioni del Calcoto Mauro Picone, Sezione di Napoli, Rapporto tecnico n. 272/03, (2005).

\bibitem{Fiorenz4}
{\sc Fiorenza A.} {\it Duality and reflexivity in grand Lebesgue
spaces.} Collect. Math. {\bf 51}(2000), 131\ - \ 148.

\bibitem{Fiorenz5}
{\sc  Fiorenza A. Karadzhov G.E. } {\it Grand and small Lebesgue
spaces and their analogs.} Consiglio Nationale Delle Ricerche,
Instituto per le Applicazioni del Calcoto Mauro Picine", Sezione
di Napoli, Rapporto tecnico 272/03(2005).

\bibitem{Iwaniec1}
{\sc  Iwaniec T. and  Sbordone C.} {\it On the integrability of the
Jacobian under minimal hypotheses.} Arch. Rat.Mech. Anal., {\bf
119}(1992), 129 \ - \ 143.

\bibitem{Iwaniec2}
{\sc Iwaniec T., Koskela P. and  Onninen J.} {\it Mapping of Finite
Distortion: Monotonicity and Continuity.} Invent. Math. {\bf
144}(2001), 507 \ - \ 531.

\bibitem{Iwaniec3}
{\sc  Iwaniec T. and  Martin G.} {\it Riesz transforms and related singular integrals.} J. Reine Angew. Math. 473
(1996), 25 \ - 5\ 7.

\bibitem{Jawerth1}
{\sc  Jawerth B. and  Milman M.} {\it Extrapolation theory with
applications.} Mem. Amer. Math. Soc. {\bf 440},  (1991).

\bibitem{Talenti1}
{\sc G.Talenti.} {\it Inequalities in Rearrangement Invariant Function Spaces. Nonlinear Analysis, Function
 Spaces and Applications.}  Prometheus, Prague, {\bf 5}, (1995), 177 \ - \ 230.


\vspace{9mm}

\bibitem{German1}
{\sc German P.} {\it The second iterate for the Navier-Stokes equation.}
arXiv:0806.4525v1 [math.AP] 27 Jun 2008

\bibitem{Ikeda1}
{\sc Masahiro Ikeda, Soichiro Katayama, and Hideaki Sunagawa.}
{\it Null structure in a system of quadratic
derivative nonlinear Schr\"odinger equations.}
arXiv:1305.3662v1 [math.AP] 16 May 2013

\bibitem{Kozachenko1}
 {\sc Kozachenko Yu. V., Ostrovsky E.I.} (1985). {\it The Banach Spaces of
      random Variables of subgaussian type.} Theory of Probab. and Math.
      Stat. (in Russian). Kiev, KSU, {\bf 32}, 43 \ - \  57.

\bibitem{Liflyand1}
{\sc Liflyand E., Ostrovsky E., Sirota L.} {\it Structural Properties of Bilateral Grand Lebesgue Spaces.}
Turk. J. Math.; {\bf 34} (2010), 207 \ - \ 219.

\bibitem{Okikiolu1}
{\sc Okikiolu G.O.} {\it Aspects of the theory of bounded Integral Operators in the $ L^p $  Spaces. }
Academic Press; London,   New Yotk; (1971).

\bibitem{Ostrovsky1}
{\sc  Ostrovsky E.I.} (1999). {\it Exponential estimations for random Fields and its
applications (in Russian).}  Moscow - Obninsk, OINPE.

\bibitem{Ostrovsky2}
{\sc  Ostrovsky E. and Sirota L.}
{\it Moment Banach spaces: theory and applications.}
HIAT Journal of Science and Engineering, {\bf C}, Volume 4, Issues 1 - 2,
pp. 233 \ -  \ 262, (2007).

\bibitem{Ostrovsky106}
 {\sc Ostrovsky E., Sirota L.}
{\it Monte-Carlo method for multiple parametric integrals
calculation and solving of linear integral Fredholm equations
of a second kind, with confidence regions in uniform norm.}
arXiv:1101.5381v1 [math.FA] 27 Jan 2011
(2007), 167-178.

\bibitem{Ostrovsky107}
 {\sc Ostrovsky E., Sirota L.} {\it Nikolskii-type inequalities for rearrangement invariant spaces.}
arXiv:0804.2311v1 [math.FA] 15 Apr 2008

\bibitem{Ostrovsky108}
{\sc Ostrovsky E., Rogover E. }
{\it Strichartz - type Inequalities for Parabolic and
Schr\"odinger Equations in rearrangement invariant Spaces.}
arXiv:0901.2715v1 [math.AP] 18 Jan 2009

\bibitem{Ostrovsky109}
 {\sc Ostrovsky E., Sirota L.} {\it Boundedness of operatore in bilateral Grand Lebesgue Spaces,
 with exact and weakly exact constant calculation.}
arXiv:1104.2963v1 [math.FA] 15 Apr 2011

\bibitem{Ostrovsky110}
{\sc Ostrovsky E., Sirota L.} {\it Multiple weight Riesz and Fourier transforms in bilateral
anisotropic Grand Lebesgue spaces.}
arXiv:1208.2392v1 [math.FA] 12 Aug 2012

\bibitem{Pichorides1}
{\sc  Pichorides S.K.} {\it On the best values of the constants in the theorems of M. Riesz, Zygmund and Kolmogorov. }
Studia Math. 44 (1972), 165 \ - \ 179.



\end{thebibliography}
\end{document}